\documentclass[10pt,letterpaper,leqno]{amsart}
\usepackage[OT1]{fontenc}
\usepackage{amsmath,amssymb,mathtools,needspace}
\usepackage[unicode,bookmarksnumbered]{hyperref}
\numberwithin{equation}{section}
\newtheorem{theorem}{Theorem}[section]
\newtheorem{proposition}[theorem]{Proposition}
\newtheorem{corollary}[theorem]{Corollary}
\newtheorem{lemma}[theorem]{Lemma}
\newtheorem{example}[theorem]{Example}
\newtheorem{conjecture}[theorem]{Conjecture}
\theoremstyle{plain}
\newtheorem{definition}[theorem]{Definition}
\theoremstyle{remark}

\newcommand{\R}{\mathbb R}

\newcommand{\eps}{\varepsilon}
\newcommand{\dd}{\,\mathrm d}
\newcommand{\Br}{\operatorname{Br}}

\DeclareMathOperator{\supp}{supp}
\DeclareMathOperator{\dist}{dist}

\hypersetup{
  pdfauthor={},
  pdftitle={An improved weighted broad L2 estimate for planar curves},
  pdfsubject={Weighted broad restriction},
  pdfkeywords={Fourier extension, broad estimates, Katz--Tao sets, Furstenberg estimates, refined decoupling}
}
\title[An improved weighted broad $L^2$ estimate]
{An improved weighted broad $L^2$ estimate in $\R^2$}
\date{}
\author{Yingyu Zhao}
\address{Yingyu Zhao, Qiuzhen College, Tsinghua University, Beijing, China}
\email{zhao-yy25@mails.tsinghua.edu.cn}
\author{Yunhao Zhao}
\address{Yunhao Zhao, School of Advanced Interdisciplinary Sciences,
University of Chinese Academy of Sciences, Beijing, China}
\email{zhaoyunhao23@mails.ucas.ac.cn}

\begin{document}
\begin{abstract}
For a compact smooth plane curve $S$ with nonzero curvature and a union $X\subseteq B_R$ of unit balls satisfying a one-dimensional Katz--Tao condition, we prove $ \|\Br_A E_S f\|_{L^2(X)} \lesssim_{\eps,\eps'}R^\eps |X|^{13/62}\|f\|_2. $ The proof relies on a refined local broad estimate that uses density decomposition and the short frequency intervals associated with coarse tubes.
\end{abstract}

\maketitle

\section{Introduction}

\subsection{Overview}\label{subsec:overview}

Let $S\subseteq\R^2$ be a compact smooth curve with nonzero curvature and $\sigma_S$ be the surface measure on $S$. The associated Fourier extension operator is
\[
 E_S f(x)=\int_S e^{ix\cdot\xi}f(\xi)\dd\sigma_S(\xi).
\]
We consider weighted $L^2$ estimates of the form
\begin{equation}\label{general-case}
    \|E_Sf\|_{L^2(X)}\lesssim |X|^\alpha\|f\|_{L^2(S)},
\end{equation}
where $X\subseteq B_R=B(0,R)$ is measurable and $\alpha>0$. By H{\"o}lder's inequality, \eqref{general-case} follows from the following weighted $(L^2,L^p)$ estimate
\begin{equation}\label{weight L^2-L^p estimate}
    \|E_Sf\|_{L^p(X)}\lesssim\|f\|_{L^2(S)},
\end{equation}
where $p$ is chosen so that $1/2-1/p=\alpha$. Conversely, if \eqref{general-case} holds for every measurable $X\subseteq\R^2$ of finite measure, taking $X$ to be the superlevel set of $|E_Sf|$ implies the weak-type version of \eqref{weight L^2-L^p estimate} with the same $\alpha$.

For general measurable set $X$, the Stein--Tomas estimate gives \eqref{general-case} with $\alpha\geqslant1/3$. To obtain a smaller $\alpha$, we can replace $E_S$ with its bilinear analog ${\rm Bil}E_Sf:=|E_{S_1}f_1E_{S_2}f_2|^{1/2}$ where $S_1,S_2\subseteq S$ are transverse arcs and $f_j=f|_{S_j}$. The bilinear restriction estimate \cite{TVV98} gives the weighted estimate with $\alpha\geqslant1/4$ if $E_S$ is replaced by ${\rm Bil}E_S$ in \eqref{general-case}. In \cite{Wu25}, Wu obtains a stronger weighted $L^2$ estimate for sets satisfying the Katz--Tao condition by replacing the extension operator $E_S$ with the broad operator. We recall the relevant definitions below.

\begin{definition}\label{def:broad}
Let $1\ll K\ll R^{1/2}$, $A\in\mathbb{N}$, and let $\Sigma:=\{\sigma\}$ be a finitely overlapping cover of $S$ by $K^{-1}$-arcs. Let $f_{\sigma}:=f\cdot\chi_{\sigma}$. For any $f:S\to\mathbb{C}$ and $x\in\mathbb{R}^2$, we define
    \begin{equation}\label{eq:broad-curve}
        {\rm Br}_{A} E_S f(x):=\max_{\substack{\Sigma'\subseteq\Sigma\\ \#\Sigma'=A}}\min_{\sigma\in\Sigma'}|E_S f_{\sigma}(x)|.
    \end{equation}
\end{definition}

\begin{definition}\label{def:KT}
Let $0<\delta<1$, $0<s\leqslant2$ and $C>0$. A finite set $E\subseteq\R^2$ is called a Katz--Tao $(\delta,s,C)$-set if
\[
 \#(E\cap B(x,r))\leqslant C(r/\delta)^s,
 \quad x\in\R^2,\ \delta\leqslant r\leqslant1.
\]
A similar definition applies to unions of $\delta$-balls by counting their centers.
\end{definition}

For a union $X\subseteq B_R$ of unit balls whose $R^{-1}$-dilate is a Katz--Tao $(R^{-1},1)$-set, Wu \cite[Theorem~1.3]{Wu25} proved
\begin{equation}\label{eq:Wu-main}
 \|\Br_A E_S f\|_{L^2(X)}
 \lesssim_{\eps,\eps'} R^\eps |X|^{2/9}\|f\|_{L^2(S)}.
\end{equation}
His proof combines refined decoupling, the dual version of the two-ends Furstenberg estimate of Wang and Wu \cite{WW24} and induction on scales.

More recently, Cao, Guo, Li and Wang \cite{CGLW26} extended Wu's weighted restriction estimates to the case where $R^{-1}X$ is a Katz--Tao $(R^{-1},s)$-set with $0\leqslant s\leqslant2$. Their proof combines refined decoupling with the two-ends Furstenberg estimates of Wang and Wu \cite{WW26}. For $s=1$, their broad estimate still has exponent $2/9$. We improve this exponent to $13/62$. More precisely, we prove the following theorem.

\begin{theorem}\label{thm:main}
Let $R\geqslant2$, $0<\eps\ll1$, and $2\leqslant K\leqslant R^{\eps^4}$. Let $X\subseteq B(0,R)$ be a union of unit balls such that $R^{-1}X$ is a Katz--Tao $(R^{-1},1)$-set. If $0<\eps'<\log K/\log R$ and $A\geqslant R^{\eps'}$, then
\begin{equation}\label{eq:main}
 \|\Br_A E_S f\|_{L^2(X)}
 \leqslant C_{\eps,\eps'}
 R^\eps |X|^{13/62}\|f\|_{L^2(S)}.
\end{equation}
\end{theorem}

For the truncated parabola $S=\{(\xi,\xi^2):\xi\in[-1,1]\}$, Wu \cite{Wu25} already gives an example showing that $\alpha\geqslant 1/6$ for the broad operator. The following example gives the same lower bound using the rectangular-lattice construction of Barcel\'o, Bennett, Carbery, Ruiz and Vilela~\cite{BBCRV07}.

\begin{example}[A lattice example]\label{ex:lattice-sixth}
There are absolute constants $c_0,C>0$ such that, for $R\gg1$ and all integers $K,A$ with $C\leqslant K\leqslant c_0R^{1/3}$ and $1\leqslant A\leqslant c_0K$, there are a function $f$ and a union $X\subseteq B_R$ of disjoint unit balls such that the $R^{-1}$-dilate of $X$ is a Katz--Tao $(R^{-1},1)$-set with a constant independent of $R$, and
\begin{equation}\label{eq:lattice-sixth}
 |X|\sim R,\quad
 \|\Br_A E_Sf\|_{L^2(X)}
 \gtrsim K^{-1}|X|^{1/6}\|f\|_{L^2(S)}.
\end{equation}
\end{example}

Let $N=\lfloor R^{1/3}\rfloor$ and fix a sufficiently small absolute constant $c>0$. Set
\[
 I_k=[k/N-c/R,k/N+c/R],\quad
 f(\xi,\xi^2)=\frac{1}{\sqrt{1+4\xi^2}}
 \sum_{k=1}^{\lfloor cN\rfloor}\mathbf{1}_{I_k}(\xi).
\]
Let $X$ be the union of unit balls centered at $P=\{(2\pi Na,2\pi N^2b)\in[-cR,cR]^2:a,b\in\mathbb Z\}$. The balls are disjoint and lie in $B_R$. The lattice spacings give $|X|\sim R$ and $\#(P\cap B(x,r))\lesssim r$ for $1\leqslant r\leqslant R$.

On the ball centered at $(2\pi Na,2\pi N^2b)$, the phase on $I_k$ is $2\pi(ak+bk^2)+O(c)$. In the stated range of $K$, there are $\gtrsim K$ arcs whose parameter intervals each contain $\gtrsim N/K$ of the $I_k$. Taking real parts after canceling the arclength density gives $|E_Sf_\sigma|\gtrsim N/(KR)$ on $X$ for these arcs. Since $\|f\|_{L^2(S)}^2\sim N/R$, for $A\leqslant c_0K$ we obtain
\[
 \|\Br_A E_Sf\|_{L^2(X)}
 \gtrsim R^{1/2}\frac{N}{KR}
 \sim K^{-1}|X|^{1/6}\|f\|_{L^2(S)}.
\]

To rule out an exponent $\alpha<1/6$ in Theorem~\ref{thm:main}, choose $\eps>0$ sufficiently small that $\alpha+\eps+\eps^4/2<1/6$ and let
\[
 K=\lfloor R^{\eps^4/2}\rfloor,\quad
 \eps'=\eps^4/4,\quad A=\lceil R^{\eps'}\rceil.
\]
For sufficiently large $R$, these parameters satisfy the hypotheses of both Theorem~\ref{thm:main} and Example~\ref{ex:lattice-sixth}. The lower bound then exceeds $R^\eps|X|^\alpha\|f\|_2$ by a fixed positive power of $R$.

We conjecture the matching upper bound.

\begin{conjecture}\label{conj:sharp-broad}
Let $S=\{(\xi,\xi^2):\xi\in[-1,1]\}$. Under the hypotheses of Theorem~\ref{thm:main}, we have
\begin{equation}\label{eq:conjectural-broad}
 \|\Br_AE_Sf\|_{L^2(X)}
 \leqslant C_{\eps,\eps'}
 R^\eps|X|^{1/6}\|f\|_{L^2(S)}.
\end{equation}
\end{conjecture}

\subsection{Outline of the proof}

Wu \cite[Section~2.3]{Wu25} obtains $2/9$ by balancing the radius bound $r^{2/5}$ with the two-ends bound $|X_1|r^{-1/2}$ for the squared broad norm. We refine the radius bound by retaining spatial density and local frequency width.

The stopping decomposition isolates the low-density part $X_\ell$, where the density parameter gives a gain in Wu's radius bound. The stopping square centers satisfy a Katz--Tao bound independent of the density threshold, so the incidence estimate also applies to the coarse tubes. On smaller stopping squares, the short frequency intervals associated with these tubes give a gain in the broad estimate.

For larger stopping squares, density layers provide two successive thresholds. The lower threshold controls the number of squares through their masses, and the upper threshold bounds the density of the layer inside each square. The resulting local estimate has the factor $r^{13/36}+|Y|^{5/12}$. The smaller radius power improves the balance with the two-ends bound, and the mass term fits within the resulting estimate. This yields $13/62$, with the one-end case handled by induction on scales.

\section{Preliminaries}\label{sec:prelim}

Throughout this paper, we use $\#E$ to denote the cardinality of a finite set. For $U,V\geqslant0$, $U\lesssim V$ means $U\leqslant CV$, and $U\sim V$ means both inequalities. We write $U\lessapprox V$ if $U\leqslant C_\eps R^\eps V$ for every $\eps>0$. The reverse inequalities are defined similarly. At local radius $r$, the same notation uses $r$ in place of $R$.

\subsection{Some reductions}

After a finite partition of the curve and a change of coordinates, we may assume that $S$ is the graph of a smooth function $\Phi:[-1,1]\to \R$ with $|\Phi''|\sim1$. Let
\begin{equation}\label{eq:extension}
 Ef(x)=\int_{-1}^1e^{i(x_1\xi+x_2\Phi(\xi))}f(\xi)\dd\xi.
\end{equation}
As in Definition~\ref{def:broad}, we consider the following equivalent definition.
\begin{definition}\label{def:broad-operator}
Let $1\ll K\ll R^{1/2}, A\in\mathbb{N}$, and let $\Sigma:=\{\sigma\}$ be the $K^{-1}$-partition of $[-1,1]$. Let $f_{\sigma}:=f\cdot\chi_{\sigma}$. For any $f:[-1,1]\to\mathbb{C}$ and $x\in\mathbb{R}^2$, we define
    \begin{equation}\label{Broad-Wu}
        {\rm Br}_{A} E f(x):=\max_{\substack{\Sigma'\subseteq\Sigma\\ \#\Sigma'=A}}\min_{\sigma\in\Sigma'}|E f_{\sigma}(x)|.
    \end{equation}
\end{definition}

The broad operator has the following standard properties.
\begin{lemma}\label{Triangle inequality for broad operator}
Let $A_1,A_2\in\mathbb{N}$ and $f_1,f_2:[-1,1]\to\mathbb{C}$. Then
    \begin{equation}\label{eq:broad-subadditivity}
        {\rm Br}_{A_1+A_2} E (f_1+f_2)\leqslant{\rm Br}_{A_1} E (f_1)+{\rm Br}_{A_2} E (f_2).
    \end{equation}
\end{lemma}
\begin{lemma}\label{pointwise estimate}

For the partition $\Sigma$ in Definition~\ref{def:broad-operator} and $f:[-1,1]\to\mathbb C$, we have
\begin{equation*}
    {\rm Br}_{A}Ef(x)\leqslant\left(\sum_{\sigma\in\Sigma}|Ef_{\sigma}(x)|^2\right)^{1/2}.
\end{equation*}
\end{lemma}

\begin{proposition}\label{main result}
Let $\varepsilon,K$ be two parameters satisfying $0<\varepsilon\ll1\ll K\leqslant R^{\varepsilon^4}\ll R$. Let $\varepsilon'\in(0,\log K/\log R)$. Let $X\subseteq B_R$ be a union of unit balls and $R^{-1}X$ be a Katz-Tao $(R^{-1},1)$ set. Then when $A\geqslant R^{\varepsilon'}$, we have
    \begin{equation}\label{weighted L2-13/62}
         \|{\rm Br}_AEf\|_{L^2(X)}\lesssim_{\varepsilon,\varepsilon'}R^{\varepsilon}|X|^{13/62}\|f\|_{L^2}.
    \end{equation}
\end{proposition}

\subsection{Wave packet decomposition}

In this section, we will first formally construct the wave packet decomposition of $Ef$ and then state some of its properties for later use. The idea of the wave packet decomposition is to break up $Ef$ into pieces $Ef_{\theta,v}$ that are "localized" in both space and frequency.

For the parameters $\eps,\eps'$ in Proposition~\ref{main result}, we take a small constant $0<\eps_0\ll\min\{\eps^{100},\eps'\}$. Cover $[-1,1]$ by finitely overlapping $R^{-1/2}$-intervals $\theta$ and denote $\Theta:=\{\theta\}$. Let $\phi_{\theta}$ be a smooth partition of unity adapted to this cover. We write $f_\theta=f\cdot\phi_{\theta}$. Next, we cover $\mathbb{R}$ by finitely overlapping $R^{1/2+\varepsilon_0}$-intervals $v$. Let $\eta_v$ be a smooth partition of unity adapted to this cover. Then $\eta_v^{\vee}$ is a Schwartz function adapted to  $[-R^{-(1/2+\varepsilon_0)}/2,R^{-(1/2+\varepsilon_0)}/2]$.

Now choose a Schwartz function $\tilde{\phi}_{\theta}$ such that $\tilde{\phi}_{\theta}=1$ on the $10^{-3}R^{-1/2}$-neighborhood of supp$\phi_{\theta}$ and is supported on ${\theta}$. We define $f_{\theta,v}:=\tilde{\phi}_{\theta}( f_{\theta}\star \eta_v^{\vee})$. By the rapid decay of $\eta_v^{\vee}$, we have the error estimate
\begin{equation*}
    \|f_{\theta,v}-f_{\theta}\star \eta_v^{\vee}\|_{L^{\infty}}\lesssim_{N} R^{-N}\|f\|_2.
\end{equation*}
Hence, we almost realize the first step
\begin{equation*}
    f=\sum_{\theta,v}f_{\theta,v}+{\text{error}}\approx \sum_{\theta,v}f_{\theta,v}.
\end{equation*}

Further, we have the wave packet decomposition as follows:
\begin{equation}\label{wpd}
    Ef(x)\approx\sum_{\theta,v}Ef_{\theta,v}(x).
\end{equation}
Each function $Ef_{\theta,v}(x)$ is called a wave packet. For each pair $(\theta,v)$, we define the spacetime tube $T_{\theta,v}:=\{(x_1,x_2)\in\mathbb{R}^2:|x_1-c_v+x_2\Phi'(c_{\theta})|\leqslant R^{1/2+\varepsilon_0},|x_2|\leqslant R\}$, where $c_v,c_{\theta}$ are the centers of $v,\theta$, respectively. Denote
\begin{equation*}
    {\rm dir}\, T_{\theta,v}:=\xi_\theta=(-\Phi'(c_{\theta}),1),\quad \mathbb{T}=\{T_{\theta,v}\},\quad f_T=f_{\theta,v}.
\end{equation*}
In each direction $\xi_\theta$, there is a family of $T_{\theta,v}\in\mathbb{T}$ tiling $B(0,R)$, denoted by $\mathbb{T}_\theta$.

\begin{lemma}
The wave packet decomposition satisfies the following properties.
    \begin{itemize}
\item[\textbf{WP1}] $\|Ef_T\|_{L^{\infty}(B(0,R)-T)}\lesssim_{N} R^{-\varepsilon_0N}\|f\|_2$.
\item[\textbf{WP2}]   $\{\xi_{\theta}\}_{\theta}$ is an $R^{-1/2}$-separated set.
\item[\textbf{WP3}] $\mathbb{T}_\theta$ is finitely overlapping.
\item[\textbf{WP4}] The functions $f_T$ are approximately orthogonal. More precisely, for any collection $\mathbb{W}\subseteq\mathbb{T}$, we have
        \begin{equation}\label{eq:packet-energy}
            \sum_{T\in\mathbb T}\|f_T\|_2^2\lesssim\|f\|_2^2,\quad
            \left\|\sum_{T\in\mathbb W}f_T\right\|_2^2\lesssim\sum_{T\in\mathbb W}\|f_T\|_2^2.
        \end{equation}
    \end{itemize}
\end{lemma}

\subsection{Main Tools}

We recall refined decoupling and the two-ends Furstenberg estimate. They address two aspects of wave packet analysis: the former controls interference among wave packets and the latter describes the incidence geometry of the associated tubes. We begin with the refined decoupling estimate of Guth, Iosevich, Ou and Wang \cite[Theorem~4.2]{GIOW20}.

\begin{theorem}\label{thm:refined}
Let $\mathbb W\subseteq\mathbb T$ and $f=\sum_{T\in\mathbb W}f_T$. Suppose that $\|f_T\|_{L^2(\mathbb{R})}$ is roughly constant for all $T\in\mathbb W$. Let $X\subseteq B_R$ be a union of $R^{1/2}$-balls each of which intersects at most $M$ tubes in $\mathbb W$. Then
\begin{equation}\label{eq:refined}
    \|Ef\|_{L^6(X)}
    \lessapprox_{\varepsilon} M^{1/3}
    \left(\sum_{T\in\mathbb W}
    \|Ef_T\|_{L^6(\omega_{B_R})}^6\right)^{1/6}.
\end{equation}
where $\omega_{B_R}(x)=(1+|x|/R)^{-100}$.
\end{theorem}

We now turn to the two-ends Furstenberg estimate of Wang and Wu \cite{WW24,WW26}. We first recall some related definitions.

\begin{definition}[Shading]
Let $T$ be a $\delta$-tube. A shading of $T$ is a set $Y(T)$ that is a union of $\delta$-balls in $T$. We refer to the pair $(\mathbb{T},Y)$ as a set of tubes and their associated shading. We say $Y$ is $\lambda$-dense if $|Y(T)|\geqslant\lambda|T|$ for all $T\in\mathbb{T}$.
\end{definition}

\begin{definition}[Two-ends condition]
Let $(\mathbb{T},Y)$ be a set of tubes and shading. Let $0<\varepsilon_2<\varepsilon_1<1$. We say $Y$ is $(\varepsilon_1,\varepsilon_2,C)$-two-ends if there exists $C>0$ such that
    \begin{equation}\label{WTE}
        |Y(T)\cap J|\leqslant C\delta^{\varepsilon_2}|Y(T)|.
    \end{equation}
The inequality holds for all $T\in\mathbb{T}$ and all $\delta\times\delta^{\varepsilon_1}$-tube segments $J\subseteq T$. When the constant $C$ is not important in the context, we say $Y$ is $(\varepsilon_1,\varepsilon_2)$ two-ends.
\end{definition}

We use point--line duality to define Katz--Tao sets of lines. Let $F$ be the one-to-one correspondence taking $(a,b)$ to the line $y=ax-b$. For a family $L$ of non-vertical lines, its dual point set is $F^{-1}(L)$. If the domain of $F$ is contained in a compact set, then $F$ is bi-Lipschitz.

\begin{definition}
We say a family of lines in $\mathbb{R}^2$ is a Katz-Tao $(\delta,s)$-set if its dual point set in $\mathbb{R}^2$ is a Katz-Tao $(\delta,s)$-set. We say a set of $\delta$-tubes $\mathbb{T}$ is a Katz-Tao $(\delta,s)$ set if the collection of their central axes is a Katz-Tao $(\delta,s)$-set.
\end{definition}

\begin{theorem}[\cite{WW26}]\label{TEF-2-thm'}
Let $(\mathbb{T},Y)$ be a set of $\delta$-separated $\delta$-tubes in $\mathbb{R}^2$ with an $(\varepsilon_1,\varepsilon_2)$-two-ends and $\lambda$-dense shading. Suppose that $\mathbb{T}$ is a Katz-Tao-$(\delta,1)$ set.  Then
    \begin{equation}\label{TEF-2'}
        \bigg|\bigcup_{T\in\mathbb{T}}Y(T)\bigg|\gtrapprox_{\varepsilon,\varepsilon_2} \delta^{\frac{\varepsilon_1}{2}}\lambda^{\frac{1}{2}}\sum_{T\in\mathbb{T}}|Y(T)|.
    \end{equation}
\end{theorem}

By the point--line duality, Theorem~\ref{TEF-2-thm'} implies the following dual version of the two-ends Furstenberg inequality.

\begin{proposition}\label{prop:incidence}
Let $K\geqslant 2$, $M\geqslant 1$, and $\varepsilon_1>0$ satisfy $\delta<K^{-1}<\delta^{\varepsilon_1}$. Let $\mathcal Q$ be a Katz--Tao $(\delta,1)$-set of $\delta$-balls in $B_1$. For each $Q\in\mathcal Q$, let $\mathcal T(Q)$ be a set of $\delta$-separated $\delta$-tubes intersecting $Q$. For a $K^{-1}$-arc $\sigma$, let $\mathcal T_\sigma(Q):=\{T\in\mathcal T(Q): \operatorname{dir}T\in\sigma\}$. Suppose that
\[
 \#\mathcal T_\sigma(Q)
 \lesssim \delta^{\varepsilon_1}\#\mathcal T(Q),
 \quad
 \#\mathcal T(Q)\geqslant M
\]
for each $Q\in\mathcal Q$ and $K^{-1}$-arc $\sigma$. Then we have
\begin{equation}\label{eq:incidence}
 \#\bigcup_{Q\in\mathcal Q}\mathcal T(Q)
 \gtrapprox_{\varepsilon,\varepsilon_1}
 K^{-1/2}M^{3/2}\delta^{1/2}\#\mathcal Q.
\end{equation}
\end{proposition}

The following lemma is \cite[Lemma~2.7]{Wu25}.

\begin{lemma}\label{lem:KT-sampling}
Let $0<\delta<1$, and let $\mathcal Q$ be a family of $\delta$-balls in $B_1$. Suppose that $\gamma_{\mathcal Q}$ is the smallest constant such that
\[
    \#\{Q\in\mathcal Q:Q\cap B(x,r)\neq\varnothing\}
    \leqslant \gamma_{\mathcal Q}\frac{r}{\delta}
\]
holds for every $x\in\mathbb R^2$ and $r\in[\delta,1]$. Then there exists a Katz--Tao $(\delta,1)$-set $\mathcal Q'\subseteq\mathcal Q$ such that $\#\mathcal Q'\gtrapprox \gamma_{\mathcal Q}^{-1}\#\mathcal Q$.
\end{lemma}

The first method in \cite[Section~2.3]{Wu25} gives the following two estimates by combining the two tools introduced above with dyadic pigeonholing steps carried out in the right order.

\begin{lemma}\label{lem:pure-radius}
There are absolute constants $\kappa_0,C>0$ such that the following holds. Let $1\leqslant r\leqslant R$, $\varepsilon'>0$, and $f\in L^2([-1,1])$. Let $X\subseteq B_r$ be a union of unit balls such that $r^{-1}X$ is a Katz--Tao $(r^{-1},1)$-set. Let $2\leqslant K\leqslant r^{\kappa_0}$ and $A\geqslant r^{\varepsilon'}$. Then
\begin{equation}\label{eq:pure-radius}
    \|\Br_A Ef\|_{L^2(X)}^2
    \lessapprox_{\varepsilon'} K^C r^{2/5}\|f\|_{L^2}^2.
\end{equation}
If $\supp f$ can be covered by $D$ intervals of length $r^{-1/2}$, then
\begin{equation}\label{eq:frequency-support}
    \|\Br_A Ef\|_{L^2(X)}^2
    \lessapprox_{\varepsilon'} K^C r^{1/4}D^{1/2}\|f\|_{L^2}^2.
\end{equation}
\end{lemma}

Under the frequency support hypothesis, the selected tube multiplicity satisfies $M\lesssim D$, so the same reduction with \cite[equation~(2.26)]{Wu25} gives \eqref{eq:frequency-support}.

\section{Low density estimate and tube counting}\label{sec:stopping}

We first decompose $X$ into a low density part and its intersections with a family of disjoint dyadic squares. We estimate the low density part using the incidence estimate and refined decoupling. We then group the wave packets into coarse tubes and combine the tube counts with the frequency support estimate to bound the sum of the local broad integrals.

In what follows, we assume $A\geqslant R^{\varepsilon'}$ and $2\leqslant K\leqslant R^\kappa$, where $\varepsilon'>0$ is fixed and $\kappa>0$ is sufficiently small.

\begin{lemma}\label{lem:density-stopping}
Let $X\subseteq B_R$ be a union of unit balls and $\gamma>0$. Let $Q_0\supseteq B_R$ be a square of side length comparable to $R$. Then there exist $X_{\ell}\subseteq X$ and a disjoint family \(\mathcal B\) of dyadic sub-squares of \(Q_0\) with \(R^{1/2}\lesssim\ell(B)\lesssim R, B\in\mathcal B\) such that
\begin{equation}\label{eq:low-density}
    X=X_\ell\sqcup\bigsqcup_{B\in\mathcal B}(X\cap B)\ \text{ and }\ |X_\ell\cap B(x,r)|\lesssim\gamma r,
 \ R^{1/2}\leqslant r\leqslant R.
\end{equation}
Writing $\mathcal B_{\rho}=\{B\in\mathcal B: \ell(B)=\rho\}$, we have
\begin{enumerate}
\item\label{eq:stop-mass} $(\#\mathcal B_\rho)\gamma\rho\leqslant |X|$;
\item\label{eq:stop-KT} $\#\{B\in\mathcal B_\rho:c_B\in B(x,r)\}\lesssim r/\rho$, $\rho\leqslant r\leqslant R$;
\item\label{eq:stop-equal-mass} $|X\cap B|\sim\gamma\rho$, $B\in\mathcal B_\rho\setminus\{Q_0\}$.
\end{enumerate}
\end{lemma}

\begin{proof}
If $|X|\geqslant\gamma\ell(Q_0)$, take $\mathcal B=\{Q_0\}$ and $X_\ell=\varnothing$. The conclusions follow immediately. Otherwise, divide $Q_0$ into its four dyadic children. At each stage, select each square $B$ satisfying $|X\cap B|\geqslant\gamma\ell(B)$, and subdivide the other squares whenever their children have side length at least $R^{1/2}$. We call these selected squares stopping squares, which form a disjoint family $\mathcal B$. Let $X_\ell=X\setminus\bigcup_{B\in\mathcal B}B$ to obtain the stated decomposition.

For $R^{1/2}\leqslant r\leqslant R$, $B(x,r)\cap Q_0$ can be covered by $O(1)$ dyadic sub-squares $Q$ of a common side length in $(r,2r]$. If $Q$ meets $X_\ell$, then the construction gives $|X\cap Q|<\gamma\ell(Q)$. Summing over the covering squares gives $|X_\ell\cap B(x,r)|\lesssim\gamma r$.

For \eqref{eq:stop-mass}, the selection criterion and disjointness give
\[
\sum_\rho(\#\mathcal B_\rho)\gamma\rho
\leqslant \sum_{B\in\mathcal B}|X\cap B|
\leqslant |X|.
\]

For \eqref{eq:stop-KT}, let $Q$ be a strict dyadic ancestor of some $B\in\mathcal B_\rho$. Then, we have
\[
\#\{B\in\mathcal B_\rho:B\subseteq Q\}\,\gamma\rho
\leqslant |X\cap Q|<\gamma\ell(Q).
\]
Cover $B(x,r)\cap Q_0$ by $O(1)$ dyadic sub-squares $Q$ with $\rho<\ell(Q)\lesssim r$, using $Q_0$ alone if $r\geqslant\ell(Q_0)$. Each $B\in\mathcal B_\rho$ with $c_B\in B(x,r)$ is contained in one of the covering squares. Summing the above estimate gives (ii).

Finally, for each nonroot $B\in\mathcal B_\rho$, its dyadic parent $B'$ has side length $2\rho$. Hence
\[
\gamma\rho\leqslant |X\cap B|
\leqslant |X\cap B'|<2\gamma\rho,
\]
which proves \eqref{eq:stop-equal-mass}.
\end{proof}

\subsection{The low density estimate}

The following argument retains the density parameter in the first method of \cite[Section~2.3]{Wu25}.

\begin{lemma}\label{lem:low}
Let $\gamma,\nu,\lambda>0$, $\mathbb W\subseteq\mathbb T$ and $f=\sum_{T\in\mathbb W}f_T$, where $\|f_T\|_2\sim\nu$ for all $T\in\mathbb W$. Let $X\subseteq B_R$ be a union of unit balls such that $\|\Br_AEf\|_{L^2(B)}\geqslant\lambda$ on each unit ball $B\subseteq X$. Suppose that there exists a union $X'\subseteq X$ of these unit balls with $|X'|\gtrsim |X|$ such that $|X'\cap B(x,r)|\lesssim\gamma r$ for all $x\in\mathbb R^2$ and $R^{1/2}\leqslant r\leqslant R.$ Then
\begin{equation}\label{eq:low}
 \lambda^2|X|\lessapprox_{\varepsilon'} K^C
 \gamma^{4/5}R^{2/5}\nu^2\#\mathbb W.
\end{equation}
\end{lemma}

\begin{proof}

\textbf{Step 1.} Dyadic pigeonholing. Choose a finitely overlapping cover $\mathcal Q$ of $X'$ by $R^{1/2}$-balls, and assign each unit ball to one covering ball containing it. For each $\sigma\in\Sigma$, let $\mathbb W_\sigma:=\{T\in\mathbb W:\operatorname{dir}T\in\sigma\}$ and $f_\sigma:=\sum_{T\in\mathbb W_\sigma}f_T$. For $Q\in\mathcal Q$ and a dyadic number $M\geqslant1$, define
\[
 \mathbb W_\sigma(Q)
 :=\{T\in\mathbb W_\sigma:T\cap2Q\neq\varnothing\},
 \
 \Sigma_M(Q)
 :=\{\sigma\in\Sigma:M\leqslant\#\mathbb W_\sigma(Q)<2M\}.
\]

By dyadic pigeonholing and Lemma~\ref{Triangle inequality for broad operator}, there exist dyadic numbers $M,H\geqslant1$, a subfamily $\mathcal Q'\subseteq\mathcal Q$, a union $X''\subseteq X'$ of retained unit balls and $A_1\gtrapprox A$ with the following properties. Writing $f_Q:=\sum_{\sigma\in\Sigma_M(Q)}f_\sigma$ for all $Q\in\mathcal Q'$, we have
\begin{enumerate}
\item $|X''\cap Q|\sim H$ for every $Q\in\mathcal Q'$;

\item $|X''|\gtrapprox |X'|$ and $|X''|\sim H\#\mathcal Q'$;

\item For every $Q\in\mathcal Q'$, $\lambda^2H \lessapprox_{\varepsilon'} \|\Br_{A_1}Ef_Q\|_{L^2(X''\cap Q)}^2$.
\end{enumerate}
Summing over $Q\in\mathcal Q'$ and using $|X'|\sim|X|$, we obtain
\[
 \lambda^2|X|
 \lessapprox_{\varepsilon'}
 \sum_{Q\in\mathcal Q'}
 \|\Br_{A_1}Ef_Q\|_{L^2(X''\cap Q)}^2.
\]

\textbf{Step 2.} The incidence estimate. For $R^{1/2}\leqslant r\leqslant R$, the low density condition gives
\[
 H\#\{Q\in\mathcal Q':Q\cap B(x,r)\neq\varnothing\}
 \lesssim |X''\cap B(x,3r)|
 \lesssim \gamma r.
\]
After dilation by $R^{-1}$, Lemma~\ref{lem:KT-sampling} yields a subfamily $\mathcal Q''\subseteq\mathcal Q'$ whose $R^{-1}$-dilate is a Katz--Tao $(R^{-1/2},1)$-set, with
\begin{equation}\label{KT-Q''}
    \#\mathcal Q''
 \gtrapprox \frac{H}{\gamma R^{1/2}}\#\mathcal Q'.
\end{equation}

For each $Q\in\mathcal Q''$, let
\[
 \mathcal T(Q)
 :=\bigcup_{\sigma\in\Sigma_M(Q)}\mathbb W_\sigma(Q).
\]
The local broad lower bound (3) in Step~1 implies $\#\Sigma_M(Q)\geqslant A_1$. Hence $\#\mathcal T(Q)\geqslant A_1M$ and for every $K^{-1}$-arc $\sigma$,
\[
 \#\mathcal T_{\sigma}(Q):=\#\{T\in\mathcal T(Q):\operatorname{dir}T\in\sigma\}
 \lesssim M
 \lesssim A_1^{-1}\#\mathcal T(Q)
 \lesssim R^{-\varepsilon'/2}\#\mathcal T(Q).
\]
Since $R^{\varepsilon'/2}\leqslant A_1\lesssim K\leqslant R^\kappa$ with $\kappa<1/2$, the scale and angular conditions in Proposition~\ref{prop:incidence} hold with $\delta=R^{-1/2}$ and sufficiently small $\varepsilon_1>0$. Applying Proposition~\ref{prop:incidence} to the $R^{-1}$-dilates of $2Q$ and absorbing the $R^{O(\varepsilon_0)}$ loss from the difference in tube widths, we obtain
\begin{equation*}
 \begin{aligned}
 \#\mathbb W
 \geqslant
 \#\bigcup_{Q\in\mathcal Q''}\mathcal T(Q)\gtrapprox_{\varepsilon'}
 K^{-C}\gamma^{-1}M^{3/2}|X|R^{-3/4},
 \end{aligned}
\end{equation*}
where the last inequality uses \eqref{KT-Q''} and $H\#\mathcal Q'\sim |X''|\gtrapprox |X|$. Consequently,
\begin{equation}\label{eq:low-incidence}
    M\lessapprox_{\varepsilon'}K^C
 \left(\gamma R^{3/4}\#\mathbb W|X|^{-1}\right)^{2/3}.
\end{equation}

\textbf{Step 3.} The analytic estimates. Cauchy--Schwarz gives $\|Ef_T\|_\infty\lesssim R^{-1/4}\nu$. For $Q\in\mathcal Q'$ and $\sigma\in\Sigma_M(Q)$, wave packet localization yields, up to an error,
\[
 \|Ef_\sigma\|_{L^\infty(Q)}
 \lesssim
 \sum_{T\in\mathbb W_\sigma(Q)}\|Ef_T\|_{L^\infty(Q)}
 \lesssim MR^{-1/4}\nu.
\]
Thus $\Br_{A_1}Ef_Q\lesssim MR^{-1/4}\nu$ on $Q$. Summing over $\mathcal Q'$ and using Step~1 and \eqref{eq:low-incidence}, we obtain
\begin{equation}\label{eq:low-triangle}
 \begin{aligned}
 \lambda^2|X|
 &\lessapprox_{\varepsilon'}
 M^2R^{-1/2}\nu^2
 \sum_{Q\in\mathcal Q'}|X''\cap Q|\lesssim |X|M^2R^{-1/2}\nu^2\\
 &\lessapprox_{\varepsilon'}K^C
 \gamma^{4/3}|X|^{-1/3}R^{1/2}
 \nu^2(\#\mathbb W)^{4/3}.
 \end{aligned}
\end{equation}

For each $\sigma\in\Sigma$, let $\mathcal Q'_\sigma :=\{Q\in\mathcal Q':\sigma\in\Sigma_M(Q)\}$. Every $Q\in\mathcal Q'_\sigma$ intersects $\lesssim M$ tubes from $\mathbb W_\sigma$. Plancherel gives $\|Ef_T\|_{L^2(\omega_{B_R})}\lesssim R^{1/2}\nu$. Interpolating with the preceding $L^\infty$ bound yields $\|Ef_T\|_{L^6(\omega_{B_R})}\lesssim\nu$. Denote $U_\sigma:=\bigcup\mathcal Q'_\sigma$. Hence Theorem~\ref{thm:refined} implies
\[
 \begin{aligned}
 \|Ef_\sigma\|_{L^6(U_\sigma)}^2
 &\lessapprox M^{2/3}
 \left(
 \sum_{T\in\mathbb W_\sigma}
 \|Ef_T\|_{L^6(\omega_{B_R})}^6
 \right)^{1/3}\lesssim M^{2/3}(\#\mathbb W_\sigma)^{1/3}\nu^2.
 \end{aligned}
\]
Using Lemma~\ref{pointwise estimate} and H\"older's inequality,
\begin{align*}
 \lambda^2|X|
 &\lessapprox_{\varepsilon'}\sum_{Q\in\mathcal Q'}
 \|\Br_{A_1}Ef_Q\|_{L^2(X''\cap Q)}^2\\
 &\lesssim\sum_{\sigma\in\Sigma}
 \|Ef_\sigma\|_{L^2(X''\cap U_\sigma)}^2\\
 &\leqslant |X''|^{2/3}\sum_{\sigma\in\Sigma}
 \|Ef_\sigma\|_{L^6(U_\sigma)}^2\\
 &\lessapprox K^C|X|^{2/3}M^{2/3}\nu^2(\#\mathbb W)^{1/3}.
\end{align*}
Applying Step~1 and \eqref{eq:low-incidence}, we obtain
\begin{equation}\label{eq:low-decoupling}
 \lambda^2|X|\lessapprox_{\varepsilon'}
 K^C\gamma^{4/9}|X|^{2/9}R^{1/3}\nu^2(\#\mathbb W)^{7/9}.
\end{equation}
Taking the geometric mean of \eqref{eq:low-triangle} and \eqref{eq:low-decoupling} with weights $2/5$ and $3/5$ gives \eqref{eq:low}.
\end{proof}

\subsection{Coarse tube counting}

Let $\mathbb W\subseteq\mathbb T$ and $f=\sum_{T\in\mathbb W}f_T$, with $\|f_T\|_2\sim\nu$, and fix $\rho\geqslant R^{1/2}$. Each wave packet is indexed by an $R^{-1/2}\times R^{1/2}$-rectangle in the direction and position parameters. Assuming $R^{O(\varepsilon_0)}\rho/R\ll K^{-1}$, we now partition this parameter space into coarser $(\rho/R)\times\rho$-rectangles, separately for each $K^{-1}$-arc $\sigma\in\Sigma$. Since $|x_2|\leqslant R$, the tubes whose parameters lie in the same $(\rho/R)\times\rho$-rectangle $\mathcal R_{\rho}$ are contained in a common $O(R^{\varepsilon_0}\rho)\times O(R)$-tube $T_{\mathcal R_{\rho}}$, which we call a coarse tube. When the rectangle need not be specified, we simply write $T_\rho$ for $T_{\mathcal R_\rho}$. Denote the family of these coarse tubes by $\mathbb T_{\rho}$, which can be partitioned into $O(1)$ families, each of which is $(\rho/R)$-separated after dilation by $R^{-1}$.

For each $ T_{\mathcal R_{\rho}}\in\mathbb T_{\rho}$, let $\mathbb W[T_{\mathcal R_{\rho}}]:=\{T_{\theta,v}\in\mathbb W:(c_\theta,c_v)\in\mathcal R_{\rho}\}$. For $F\geqslant1$, write $\mathbb T_{\rho,F}:=\{T_{\rho}\in\mathbb T_{\rho}:\#\mathbb W[T_{\rho}]\sim F\}$. Then we have $F\#\mathbb T_{\rho,F}\lesssim\#\mathbb W$. Fix $m\geqslant1$ and $\mathbb V\subseteq\mathbb T_{\rho,F}$. For each $\rho$-square $B\subseteq Q_0$, define
\[
\begin{aligned}
\mathbb V_\sigma(B)
&:=
\{T_{\rho}\in\mathbb V:
T_{\rho}\cap B\neq\varnothing,\ 
\operatorname{dir}T_{\rho}\in\sigma\},\\
\Sigma^+_m(B)
&:=
\{\sigma\in\Sigma:\#\mathbb V_\sigma(B)\geqslant m\}.
\end{aligned}
\]
For each such $B$ and each $\sigma\in\Sigma^+_m(B)$, choose $\mathbb V'_\sigma(B)\subseteq\mathbb V_\sigma(B)$ with $\#\mathbb V'_\sigma(B)\sim m$.

\begin{lemma}\label{lem:coarse-tubes}
Suppose that $\#\Sigma^+_m(B)\geqslant R^{\varepsilon'} $ for every $B\in\mathcal B_\rho$. Then,
\begin{equation}\label{eq:coarse-counts}
\begin{aligned}
\#\mathbb V
\gtrapprox_{\varepsilon'}
K^{-C}m(\#\mathcal B_\rho)^{1/2},\quad
\#\mathbb V
\gtrapprox_{\varepsilon'}
K^{-C}m^{3/2}(\rho/R)^{1/2}\#\mathcal B_\rho.
\end{aligned}
\end{equation}
\end{lemma}

\begin{proof}
For every $B\in\mathcal B_\rho$, the selected families satisfy
\[
\sum_{\sigma\in\Sigma_m^+(B)}
\#\mathbb V'_\sigma(B)
\sim m\,\#\Sigma_m^+(B)
\gtrsim mR^{\varepsilon'}.
\]
\textbf{Step 1.} Double counting. For each $B\in\mathcal B_\rho$, choose $\sigma_1,\sigma_2\in\Sigma_m^+(B)$ with $\operatorname{dist}(\sigma_1,\sigma_2)\gtrsim K^{-1}$. Then every pair $T_{\rho,1}\in\mathbb V'_{\sigma_1}(B),\,T_{\rho,2}\in\mathbb V'_{\sigma_2}(B)$ meets $B$ and satisfies $\angle(T_{\rho,1},T_{\rho,2})\gtrsim K^{-1}$. Hence
\[
\sum_{B\in\mathcal B_\rho}
\sum_{\substack{
T_{\rho,1}\in\mathbb V'_{\sigma_1}(B)
}}\sum_{T_{\rho,2}\in\mathbb V'_{\sigma_2}(B)}
1
\gtrsim m^2\#\mathcal B_\rho.
\]

Conversely, fix $T_{\rho,1},T_{\rho,2}\in\mathbb V$ with $\angle(T_{\rho,1},T_{\rho,2})\gtrsim K^{-1}$. If both tubes meet a square $B\in\mathcal B_\rho$, then $c_B$ lies in the intersection of their $O(\rho)$-neighborhoods. Since each tube has width $O(R^{\varepsilon_0}\rho)$, this intersection is contained in a ball of radius $O(KR^{\varepsilon_0}\rho)$. Therefore, the disjointness of the $\rho$-squares gives
\[
\#\{B\in\mathcal B_\rho:
T_{\rho,1}\cap B\neq\varnothing,\
T_{\rho,2}\cap B\neq\varnothing\}
\lesssim K^2R^{2\varepsilon_0}.
\]
Double counting yields
\[
m^2\#\mathcal B_\rho
\lesssim K^2R^{2\varepsilon_0}(\#\mathbb V)^2,
\]
which proves the first bound in \eqref{eq:coarse-counts}.

\textbf{Step 2.} Two-ends Furstenberg. Let $\delta=\rho/R$ and dilate by $R^{-1}$. For each $B\in\mathcal B_\rho$, let $\widehat B$ be a ball centered at $c_B/R$ of radius comparable to $\delta$ containing $R^{-1}B$. Write $\widehat{\mathcal B}_\delta :=\{\widehat B:B\in\mathcal B_\rho\}$. By \eqref{eq:stop-KT},
\[
\#\{\widehat B\in\widehat{\mathcal B}_\delta:
c_{\widehat B}\in B(y,s)\}
\lesssim \frac{s}{\delta},
\quad \delta\leqslant s\leqslant1.
\]
Thus, $\widehat{\mathcal B}_\delta$ is a Katz--Tao $(\delta,1)$-set, with $\#\widehat{\mathcal B}_\delta=\#\mathcal B_\rho$.

Discretizing the rescaled coarse tubes at scale $\delta$ costs $R^{O(\varepsilon_0)}$. Let $\mathcal T(\widehat B)$ be the resulting $\delta$-separated family of $\delta$-tubes meeting $\widehat B$ from the selected families $\mathbb V'_\sigma(B)$, $\sigma\in\Sigma_m^+(B)$. Since $\delta\ll K^{-1}$, taking $\varepsilon_0$ sufficiently small compared with $\varepsilon'$ gives
\[
\#\mathcal T(\widehat B)\gtrsim mR^{\varepsilon'/2},
\quad
\#\mathcal T_\sigma(\widehat B)
\lesssim m
\lesssim R^{-\varepsilon'/2}\#\mathcal T(\widehat B)
\]
for every $K^{-1}$-arc $\sigma$. Take $\varepsilon_1=\varepsilon'/2$. By hypothesis, $R^{\varepsilon'}\leqslant \#\Sigma_m^+(B)\leqslant\#\Sigma\lesssim K$. Together with the assumption $R^{O(\varepsilon_0)}\delta\ll K^{-1}$, this gives $\delta<K^{-1}<\delta^{\varepsilon_1}$ for sufficiently large $R$. Each coarse tube gives rise to at most $R^{O(\varepsilon_0)}$ covering tubes. Applying Proposition~\ref{prop:incidence}, we obtain
\[
\begin{aligned}
R^{O(\varepsilon_0)}\#\mathbb V
&\gtrsim
\#\bigcup_{\widehat B\in\widehat{\mathcal B}_\delta}
\mathcal T(\widehat B)\gtrapprox_{\varepsilon'}
K^{-1/2}
\bigl(mR^{\varepsilon'/2}\bigr)^{3/2}
\delta^{1/2}\#\widehat{\mathcal B}_\delta.
\end{aligned}
\]
Since $\varepsilon_0$ is sufficiently small compared with $\varepsilon'$, the factor $R^{3\varepsilon'/4}$ dominates the preceding $R^{O(\varepsilon_0)}$ loss. Hence,
\[
\#\mathbb V
\gtrapprox_{\varepsilon'}
K^{-C}m^{3/2}
\left(\frac{\rho}{R}\right)^{1/2}
\#\mathcal B_\rho.
\]
\end{proof}

For each $T_{\rho}\in\mathbb T_{\rho,F}$ and $B\in\mathcal B_\rho$, we define
\[
f_{T_{\rho}}:=
\sum_{T\in\mathbb W[T_{\rho}]}f_T,\quad f_B:=
\sum_{\sigma\in\Sigma^+_m(B)}
\sum_{T_{\rho}\in\mathbb V'_\sigma(B)}
f_{T_{\rho}}.
\]
Since $\#\Sigma\lesssim K$, \eqref{eq:packet-energy} gives
\begin{equation}\label{eq:coarse-local-data}
\|f_B\|_2^2
\lessapprox_{\varepsilon'}K^CFm\nu^2.
\end{equation}
Each $f_{T_{\rho}}$ is supported in an interval of length $\sim\rho/R$. Since $f_B$ is a sum of $O(Km)$ such functions, its support can be covered by $\lessapprox_{\varepsilon'}K^Cm$ intervals of length $\sim\rho/R$.

\begin{corollary}\label{cor:coarse-frequency-final}
For each $B\in\mathcal B_\rho$, suppose that $\#\Sigma_m^+(B)\geqslant R^{\varepsilon'}$ and let $Y_B\subseteq B$ be a union of unit balls such that $\rho^{-1}Y_B$ is a Katz--Tao $(\rho^{-1},1)$-set. Then
\begin{equation}\label{eq:coarse-frequency-final}
 H(\rho):=\sum_{B\in\mathcal B_\rho}\|\Br_A Ef_B\|_{L^2(Y_B)}^2
 \lessapprox_{\varepsilon'} K^C
 \left(R^{1/2}\rho^{-1/4}+\rho^{1/2}\right)F\nu^2\#\mathbb V.
\end{equation}
\end{corollary}

\begin{proof}
The support of $f_B$ can be covered by $\lessapprox_{\varepsilon'}K^Cm$ intervals of length $O(\rho/R)$ and each $\rho/R$-interval can be covered by $O(1+\rho^{3/2}/R)$ intervals of length $\rho^{-1/2}$. The density assumption on $Y_B$ shows that $\rho^{-1}Y_B$ is a Katz--Tao $(\rho^{-1},1)$-set. Hence, by Lemma~\ref{lem:pure-radius} and \eqref{eq:coarse-local-data} we obtain
\[
\begin{aligned}
\|\Br_A Ef_B\|_{L^2(Y_B)}^2
&\lessapprox_{\varepsilon'}
K^C\rho^{1/4}m^{1/2}
\left(1+\rho^{3/2}R^{-1}\right)^{1/2}
\|f_B\|_2^2\\
&\lessapprox_{\varepsilon'}
K^C\left(\rho^{1/4}+\rho R^{-1/2}\right)
m^{3/2}F\nu^2.
\end{aligned}
\]
Summing over $B\in\mathcal B_\rho$ and using the second bound in \eqref{eq:coarse-counts}, we obtain
\[
\begin{aligned}
H(\rho)\lessapprox_{\varepsilon'}
K^C\left(R^{1/2}\rho^{-1/4}+\rho^{1/2}\right)
F\nu^2\#\mathbb V,
\end{aligned}
\]
\end{proof}

\section{A local mass estimate}\label{sec:mixed}

We next estimate $H(\rho)$ in terms of the masses $|Y_B|$. Together with the coarse tube counts and the density bounds on larger stopping squares, this will yield Proposition~\ref{prop:local-mass}.

\subsection{Estimates on stopping squares}

We begin with the following weighted local constancy estimate.

\begin{lemma}\label{lem:weighted-local-constancy}
Let $0<w\leqslant1$. Let $Y\subseteq B_r$ be a union of unit balls such that $r^{-1}Y$ is a Katz--Tao $(r^{-1},1)$-set. If $f\in L^2(\R^2)$ satisfies $\supp\widehat f\subseteq B(0,Cw)$, then
\begin{equation}\label{eq:weighted-local-constancy}
 \int_Y|f|^2\lesssim_C w\|f\|_{L^2(\R^2)}^2.
\end{equation}
\end{lemma}

\begin{proof}
Choose a Schwartz function $\phi$ such that $\widehat\phi=1$ on $B(0,C)$. Since $f=f*(w^2\phi(w\cdot))$, Cauchy--Schwarz gives
\begin{equation}\label{f-local-constancy}
    |f(x)|^2\lesssim_C
 \int_{\mathbb R^2}
 w^2(1+w|x-y|)^{-10}|f(y)|^2\,dy.
\end{equation}
Since $w^{-1}\geqslant1$, the density bound on $Y$ yields
\[
 \int_Y w^2(1+w|x-y|)^{-10}\,dx
 \lesssim
 w^2\sum_{j\geqslant0}2^{-10j}
 |Y\cap B(y,2^{j+1}/w)|
 \lesssim w.
\]
Integrating \eqref{f-local-constancy} over $Y$ and applying Fubini's theorem proves \eqref{eq:weighted-local-constancy}.
\end{proof}

\begin{lemma}\label{lem:thin-band}
Let $d\geqslant1$, $0<w\leqslant1$ and $Y$ be as in Lemma~\ref{lem:weighted-local-constancy}. Suppose that $f\in L^2([-1,1])$ is supported in a union of $d$ intervals of length at most $w$. If $A$ is sufficiently large, then
\begin{equation}\label{eq:thin-band}
 \|\Br_A Ef\|_{L^2(Y)}^2
 \lesssim K^Cd|Y|^{1/2}w^{1/2}\|f\|_2^2.
\end{equation}
\end{lemma}

\begin{proof}
For sufficiently large $A$, at each $x$ we can choose $\sigma,\tau\in\Sigma$ with $\dist(\sigma,\tau)\gtrsim K^{-1}$ such that $|\Br_A Ef(x)|^2\leqslant|Ef_\sigma(x)Ef_\tau(x)|$. Thus,
\begin{equation}\label{sum1}
    \|\Br_A Ef\|_{L^2(Y)}^2
 \leqslant
 \sum_{\dist(\sigma,\tau)\gtrsim K^{-1}}
 \|Ef_\sigma Ef_\tau\|_{L^1(Y)}.
\end{equation}

Fix such a pair $\sigma,\tau$. Partition $[-1,1]$ into disjoint intervals $I_i$ of length $w$ and let $f_{\sigma,i}:=f_\sigma\mathbf 1_{I_i}, \, f_{\tau,i}:=f_\tau\mathbf 1_{I_i}.$ Each decomposition has $O(d)$ nonzero terms, and disjointness gives
\[
 \sum_i\|f_{\sigma,i}\|_2^2=\|f_\sigma\|_2^2,
 \quad
 \sum_j\|f_{\tau,j}\|_2^2=\|f_\tau\|_2^2.
\]
Let $h_{ij}:=Ef_{\sigma,i}Ef_{\tau,j}$. We have
\[
 h_{ij}(x)
 =
 \iint_{\sigma\times\tau}
 e^{i[x_1(\xi+\eta)+x_2(\Phi(\xi)+\Phi(\eta))]}
 f_{\sigma,i}(\xi)f_{\tau,j}(\eta)\,d\xi\,d\eta.
\]
The change of variables $u_1=\xi+\eta,\, u_2=\Phi(\xi)+\Phi(\eta)$ is injective on $\sigma\times\tau$, with Jacobian $|\Phi'(\eta)-\Phi'(\xi)|\gtrsim|\eta-\xi|\gtrsim K^{-1}$. Hence,
\[
 h_{ij}(x)=\int_{\mathbb R^2}e^{ix\cdot u}G_{ij}(u)\,du,
 \quad
 G_{ij}(u)=
 \frac{f_{\sigma,i}(\xi(u))f_{\tau,j}(\eta(u))}
 {|\Phi'(\eta(u))-\Phi'(\xi(u))|},
\]
where $G_{ij}$ is extended by zero outside the image of $\sigma\times\tau$. Here $\supp\widehat h_{ij}=\supp G_{ij}$ lies in a ball of radius $O(w)$. Plancherel gives
\begin{equation}\label{eq:band-bilinear}
 \begin{aligned}
 \|h_{ij}\|_{L^2(\mathbb R^2)}^2
 &\sim\|G_{ij}\|_2^2
 \lesssim K\|f_{\sigma,i}\|_2^2\|f_{\tau,j}\|_2^2.
 \end{aligned}
\end{equation}
Lemma~\ref{lem:weighted-local-constancy} and Cauchy--Schwarz give
\[
 \begin{aligned}
 \|h_{ij}\|_{L^1(Y)}
  \lesssim K^{1/2}|Y|^{1/2}w^{1/2}
 \|f_{\sigma,i}\|_2\|f_{\tau,j}\|_2.
 \end{aligned}
\]
Summing over $i,j$ and applying Cauchy--Schwarz, we obtain
\begin{equation}\label{eq:band-pair-L1}
 \begin{aligned}
 \|Ef_\sigma Ef_\tau\|_{L^1(Y)}
 \lesssim K^{1/2}d|Y|^{1/2}w^{1/2}
 \|f_\sigma\|_2\|f_\tau\|_2.
 \end{aligned}
\end{equation}
Plugging \eqref{eq:band-pair-L1} into \eqref{sum1} and using Cauchy--Schwarz again, we obtain
\[
 \begin{aligned}
 \|\Br_A Ef\|_{L^2(Y)}^2\lesssim K^{3/2}d|Y|^{1/2}w^{1/2}\|f\|_2^2,
 \end{aligned}
\]
\end{proof}

\begin{corollary}\label{cor:thin-small}
Assume the hypotheses of Corollary~\ref{cor:coarse-frequency-final}. If $R^{1/2}\leqslant \rho\leqslant R^{2/3}$ and $|Y_B|\lesssim\gamma\rho$ for each $B\in\mathcal B_\rho$, then
\begin{equation}\label{eq:thin-small}
 H(\rho)\lessapprox_{\varepsilon'} K^C
 \gamma^{5/36}R^{13/36}F\nu^2\#\mathbb V.
\end{equation}
\end{corollary}

\begin{proof}
Since $\supp f_B$ can be covered by $O(Km)$ intervals of length $O(\rho/R)$ for all $B\in\mathcal B_\rho$, by Lemma~\ref{lem:thin-band} and \eqref{eq:coarse-local-data}, we have
\begin{align*}
 \|\Br_A Ef_B\|_{L^2(Y_B)}^2
 &\lesssim K^C m |Y_B|^{1/2}(\rho/R)^{1/2}\|f_B\|_2^2\\
 &\lessapprox_{\varepsilon'} K^C
 \gamma^{1/2}m^2\rho R^{-1/2}F\nu^2.
\end{align*}
Summing over $B$ and using the second inequality in \eqref{eq:coarse-counts}, we obtain
\begin{equation}\label{eq:thin-high1}
 \begin{aligned}
 H(\rho)
 &\lessapprox_{\varepsilon'} K^C\gamma^{1/2}m^{1/2}\rho^{1/2}F\nu^2\#\mathbb V.
 \end{aligned}
\end{equation}

On the other hand, using Lemma~\ref{lem:pure-radius} at radius $\rho$ and \eqref{eq:coarse-local-data} gives
\[
 \|\Br_A Ef_B\|_{L^2(Y_B)}^2
 \lessapprox_{\varepsilon'} K^C\rho^{2/5}Fm\nu^2.
\]
Summing over $B$ and using the second inequality in \eqref{eq:coarse-counts} again, we obtain
\begin{equation}\label{thin-high2}
    H(\rho)
 \lessapprox_{\varepsilon'} K^C R^{1/2}m^{-1/2}\rho^{-1/10}F\nu^2\#\mathbb V.
\end{equation}
Taking the geometric mean of \eqref{eq:thin-high1} and \eqref{thin-high2} gives
\begin{equation}\label{quan3}
    H(\rho)
 \lessapprox_{\varepsilon'} K^C\gamma^{1/4}R^{1/4}\rho^{1/5}F\nu^2\#\mathbb V.
\end{equation}
Since $\rho\leqslant R^{2/3}$, \eqref{eq:coarse-frequency-final} gives
\begin{equation}\label{quan4}
    H(\rho)\lessapprox_{\varepsilon'} K^C R^{1/2}\rho^{-1/4}F\nu^2\#\mathbb V.
\end{equation}
Taking the geometric mean of \eqref{quan3} and \eqref{quan4} with weights $5/9$ and $4/9$, respectively, we obtain \eqref{eq:thin-small}.
\end{proof}

We now treat the larger stopping squares.

\begin{lemma}\label{lem:middle-thirteen-final}
Let $R^{2/3}\leqslant\rho\leqslant R$ and let $Y\subseteq B_\rho$ be a union of unit balls such that $\rho^{-1}Y$ is a Katz--Tao $(\rho^{-1},1)$-set. Suppose that for some $\gamma>0$,
\begin{equation}\label{eq:middle-density-cap-final}
 |Y\cap B(x,s)|\leqslant\gamma s,
 \quad x\in\mathbb R^2,\ s\geqslant R^{1/2}.
\end{equation}
Then for $f\in L^2([-1,1])$ we have
\begin{equation}\label{eq:middle-thirteen-final}
 \|\Br_A Ef\|_{L^2(Y)}^2
 \lessapprox_{\varepsilon'} K^C
 \left(\gamma^{4/5}\rho^{2/5}+\gamma^{5/36}\rho^{13/36}+R^{1/4}\right)\|f\|_2^2.
\end{equation}
\end{lemma}

\begin{proof}
Decompose each $f_\sigma$ into wave packets at scale $\rho$. By dyadic pigeonholing and Lemma~\ref{Triangle inequality for broad operator}, we may assume that there exist $\nu,\lambda>0$ and a collection of $\rho^{1/2}\times\rho$-tubes $\mathbb W_{\rho}$ such that
\[
 f=\sum_{T_\rho\in\mathbb W_\rho}f_{T_\rho},\quad
 \|f_{T_\rho}\|_2\sim\nu,\ T_\rho\in\mathbb W_\rho,
\]
and $\|\Br_A Ef\|_{L^2(B)}\sim\lambda$ for each unit ball $B\subseteq Y$.

Using four shifted dyadic grids, we may assume that every unit ball of $Y$ lies wholly in a square at the smallest stopping scale. Apply Lemma~\ref{lem:density-stopping} to $Y$ at radius $\rho$ with parameter $2\gamma$. Write $\mathcal B'_\tau$ for the stopping squares of side length $\tau$ and let $Y_B=Y\cap B$, $B\in\mathcal B'$. Then $Y=Y_\ell\sqcup\bigsqcup_\tau\bigsqcup_{B\in\mathcal B'_\tau}Y_B$ with the following three properties:
\begin{enumerate}
\item $|Y_\ell\cap B(x,r)| \lesssim \gamma r, \quad\rho^{1/2}\leqslant r\leqslant\rho$;

\item $\#\{B\in\mathcal B'_\tau:c_B\in B(x,r)\} \lesssim r/\tau, \quad\tau\leqslant r\leqslant\rho$;

\item $|Y_B| \lesssim \gamma\tau, \quad B\in\mathcal B'_\tau$.
\end{enumerate}
By \eqref{eq:middle-density-cap-final}, any square $B$ of side length $\tau\geqslant R^{1/2}$ has $|Y\cap B|\leqslant\gamma\tau<2\gamma\tau$. Thus every stopping side satisfies $\rho^{1/2}\leqslant\tau<R^{1/2}\leqslant\rho^{3/4}$. Applying Lemma~\ref{lem:low} with $X=X'=Y_\ell$ gives
\begin{equation}\label{o-quan1}
    \|\Br_A Ef\|_{L^2(Y_\ell)}^2
 \lesssim \lambda^2|Y_\ell|
 \lessapprox_{\varepsilon'} K^C\gamma^{4/5}\rho^{2/5}
 \nu^2\#\mathbb W_\rho.
\end{equation}

For each stopping side $\tau$, the bounds $K\leqslant R^\kappa \leqslant\rho^{3\kappa/2}$ and $\tau\leqslant\rho^{3/4}$ give
\[
 K\rho^{O(\varepsilon_0)}\frac{\tau}{\rho}
 \lesssim \rho^{-1/4+3\kappa/2+O(\varepsilon_0)}\ll1.
\]
We use the coarse tube construction of Section~\ref{sec:stopping} with ambient radius $\rho$ and stopping side $\tau$. For fixed $B\in\mathcal B'_\tau$ and dyadic numbers $F,m$, we define $\mathbb T_{\tau,F,\sigma}(B):=\{T_\tau\in\mathbb T_{\tau,F}:T_\tau\cap B\neq\varnothing,{\rm dir}T_\tau\in\sigma\}$ and
\[
\begin{aligned}
 \Sigma_{m,F}(B)
 &:=\{\sigma\in\Sigma:m\leqslant\#\mathbb T_{\tau,F,\sigma}(B)<2m\},\\
 f_{B,F,m}
 &:=\sum_{\sigma\in\Sigma_{m,F}(B)}
      \sum_{T_\tau\in\mathbb T_{\tau,F,\sigma}(B)}f_{T_\tau}.
\end{aligned}
\]
Then, we have
\[
 \sum_{\substack{T_\tau\in\mathbb T_\tau\\T_\tau\cap B\neq\varnothing}}
 f_{T_\tau}
 =\sum_{F,m}f_{B,F,m}.
\]
Applying Lemma~\ref{Triangle inequality for broad operator} and dyadic pigeonholing in $F,m$, we obtain a pair $F,m$ and $A_1\gtrsim A/(\log R)^C$ such that
\[
 \|\Br_A Ef\|_{L^2(Y_B)}^2
 \lessapprox_{\varepsilon'}
 \|\Br_{A_1}E f_{B,F,m}\|_{L^2(Y_B)}^2.
\]
For this choice of $F$, denote $\mathbb V=\mathbb T_{\tau,F}$ and $\Sigma_m(B)=\Sigma_{m,F}(B)$. The right-hand side vanishes if $\#\Sigma_m(B)<A_1$, so we may restrict to squares satisfying
\[
\#\Sigma_m^+(B)\geqslant\#\Sigma_m(B)\geqslant A_1
\geqslant R^{\varepsilon'/2}\geqslant\rho^{\varepsilon'/2}.
\]
Take $\mathbb V'_\sigma(B)=\mathbb V_\sigma(B)$ for $\sigma\in\Sigma_m(B)$. On the remaining arcs in $\Sigma_m^+(B)$, let $\mathbb V'_\sigma(B)$ consist of any $m$ tubes from $\mathbb V_\sigma(B)$. Define
\begin{equation*}
    f_B:=\sum_{\sigma\in\Sigma_m^+(B)}\sum_{T_\tau\in\mathbb V'_\sigma(B)} f_{T_\tau}.
\end{equation*}
Then, we have the pointwise bound $\Br_{A_1}Ef_{B,F,m}\leqslant \Br_{A_1}Ef_B$. Denote
\begin{equation*}
    H(\tau):=\sum_{B\in\mathcal B_\tau'}\|\Br_{A_1}Ef_B\|_{L^2(Y_B)}^2.
\end{equation*}

Each $\tau^{-1}Y_B$ is a Katz--Tao $(\tau^{-1},1)$-set by the condition on $\rho^{-1}Y$. Since $\tau\geqslant R^{1/3}$, we have $K\leqslant\tau^{3\kappa}$. For $\rho^{2/3}<\tau<R^{1/2}$, Corollary~\ref{cor:coarse-frequency-final} yields
\begin{equation*}
    H(\tau)\lessapprox_{\varepsilon'}K^C
 \bigl(\rho^{1/2}\tau^{-1/4}+\tau^{1/2}\bigr)F\nu^2\#\mathbb V
 \lesssim K^CR^{1/4}F\nu^2\#\mathbb V.
\end{equation*}
On $\rho^{1/2}\leqslant\tau\leqslant\rho^{2/3}$, we have the additional condition $|Y_B|\lesssim\gamma\tau$ which allows us to apply Corollary~\ref{cor:thin-small}:
\begin{equation*}
    H(\tau)\lessapprox_{\varepsilon'}K^C
 \gamma^{5/36}\rho^{13/36}F\nu^2\#\mathbb V.
\end{equation*}
Using $F\#\mathbb V\leqslant\#\mathbb W_\rho$, we sum over all dyadic stopping squares to obtain
\begin{equation}\label{o-quan2}
    \|\Br_A Ef\|_{L^2(Y\setminus Y_\ell)}^2
 \lessapprox_{\varepsilon'}\sum_{\tau,F,m} H(\tau)
 \lessapprox_{\varepsilon'} K^C
 \bigl(\gamma^{5/36}\rho^{13/36}+R^{1/4}\bigr)\nu^2\#\mathbb W_\rho.
\end{equation}
By \eqref{eq:packet-energy}, we have $\nu^2\#\mathbb W_{\rho}\lesssim\|f\|_2^2$. Combining \eqref{o-quan1} and \eqref{o-quan2}, we obtain \eqref{eq:middle-thirteen-final}.
\end{proof}

\subsection{Density layers}\label{sec:layered-mass-final}

We combine the preceding stopping-square estimates by decomposing the set into density layers.

\begin{proposition}\label{prop:local-mass}
Let $f\in L^2([-1,1])$ and $Y\subseteq B_r$ be a union of unit balls such that $r^{-1}Y$ is a Katz--Tao $(r^{-1},1,1)$-set. If $2\leqslant K\leqslant r^\kappa$ and $A\geqslant r^{\varepsilon'}$, then for every $\eps>0$,
\begin{equation}\label{eq:floor-thirteen-mass}
 \|\Br_A Ef\|_{L^2(Y)}^2
 \lesssim_{\eps,\eps'} K^C r^{\eps}
 \left(r^{13/36}+|Y|^{5/12}\right)\|f\|_2^2.
\end{equation}
\end{proposition}

\begin{proof}
Decompose each $f_\sigma$ into wave packets at scale $r$. Using the translated grids and dyadic pigeonholing as in the proof of Lemma~\ref{lem:middle-thirteen-final}, we may assume that $f=\sum_{T\in\mathbb W}f_T$ with $\|f_T\|_2\sim\nu$ and $\|\Br_AEf\|_{L^2(B)}\sim\lambda$ for each unit ball $B\subseteq Y$.

\textbf{Step 1.} Density layers. Fix a small $\theta>0$ and choose
\[
 r^{-2}\sim\gamma_0<\gamma_1<\cdots<\gamma_J,\quad
 \gamma_{j+1}/\gamma_j\leqslant r^\theta,\quad
 J=O_\theta(1),
\]
with $\gamma_J\geqslant2$. Apply Lemma~\ref{lem:density-stopping} to $Y$ at each threshold $\gamma_j$. Since the root $Q_0$ stops at $\gamma_0$ and the Katz--Tao condition rules out stopping at $\gamma_J$, the low density sets $Y_{\ell,j}$ increase from $\varnothing$ to $Y$. Thus,
\begin{equation}\label{eq:finite-layers-final}
 Y=\bigsqcup_{j=0}^{J-1}Y_j,\quad
 Y_j:=Y_{\ell,j+1}\setminus Y_{\ell,j}.
\end{equation}
Fix a nonempty $Y_j$. Since $Y_j\subseteq Y_{\ell,j+1}\subseteq B_r$, \eqref{eq:low-density} gives
\begin{equation}\label{eq:layer-density}
 |Y_j\cap B(x,t)|\lesssim \gamma_{j+1}t,
 \quad x\in\mathbb R^2,\ t\geqslant r^{1/2}.
\end{equation}
Applying Lemma~\ref{lem:low} with $X=X'=Y_j$ gives
\begin{equation}\label{eq:layer-low-final}
 \|\Br_A Ef\|_{L^2(Y_j)}^2
 \lessapprox_{\varepsilon'} K^C\gamma_{j+1}^{4/5}r^{2/5}\nu^2\#\mathbb W.
\end{equation}

\textbf{Step 2.} Estimates on stopping squares. Let $Y_B:=Y_j\cap B$ for the $\gamma_j$-stopping squares. If a stopping square of side length $\rho$ is the root or satisfies $\rho>r^{1-C_2\varepsilon_0}K^{-C_1}$, the stopping condition gives $\gamma_j\leqslant |Y|/\rho\leqslant K^{C_1}r^{C_2\varepsilon_0}|Y|/r$. Combining this with \eqref{eq:layer-low-final}, we obtain
\begin{equation}\label{eq:layer-root-final}
 \|\Br_A Ef\|_{L^2(Y_j)}^2
 \lessapprox_{\varepsilon'} K^C
 \left(\frac{\gamma_{j+1}}{\gamma_j}\right)^{4/5}
 |Y|^{4/5}r^{-2/5}\nu^2\#\mathbb W.
\end{equation}

Suppose now that all stopping sides satisfy $\rho\leqslant r^{1-C_2\varepsilon_0}K^{-C_1}$. For $C_1,C_2$ sufficiently large, this gives $r^{O(\varepsilon_0)}\rho/r\ll K^{-1}$. Apply the coarse tube selections from the proof of Lemma~\ref{lem:middle-thirteen-final} to the sets $Y_B$. For each fixed $\rho,F,m$, denote the retained $\gamma_j$-stopping squares of side length $\rho$ by $\mathcal B_\rho$ and the corresponding coarse tube family by $\mathbb V$. Then
\[
 \|\Br_A Ef\|_{L^2(Y_j)}^2
 \lessapprox_{\varepsilon'} \sum_{\rho,F,m}H(\rho).
\]
Each $\rho^{-1}Y_B$ is a Katz--Tao $(\rho^{-1},1)$-set by the assumption on $Y$. For $r^{1/2}\leqslant\rho\leqslant r^{2/3}$, Lemma~\ref{lem:density-stopping}\eqref{eq:stop-equal-mass} gives $|Y_B|\lesssim\gamma_j\rho$. By Corollary~\ref{cor:thin-small},
\begin{equation}\label{4.20}
    H(\rho)\lessapprox_{\varepsilon'} K^C
 \gamma_j^{5/36}r^{13/36}F\nu^2\#\mathbb V
 \lesssim K^C r^{13/36}F\nu^2\#\mathbb V.
\end{equation}

For $\rho>r^{2/3}$, $Y_B\subseteq Y_j$ inherits \eqref{eq:layer-density}. Apply Lemma~\ref{lem:middle-thirteen-final} to $f_B$ on $Y_B\subseteq B(c_B,\rho)$. Summing over $\mathcal B_\rho$ and using \eqref{eq:coarse-local-data}, the first inequality in Lemma~\ref{lem:coarse-tubes} and Lemma~\ref{lem:density-stopping}\eqref{eq:stop-mass}, we have
\[
 \begin{aligned}
 H(\rho)
 &\lessapprox_{\varepsilon'} K^C
 \left(\gamma_{j+1}^{4/5}\rho^{2/5}
 +\gamma_{j+1}^{5/36}\rho^{13/36}+r^{1/4}\right)Fm\nu^2\#\mathcal B_\rho\\
 &\lessapprox_{\varepsilon'} K^C
 \left(\frac{|Y|}{\gamma_j}\right)^{1/2}
 \left(\gamma_{j+1}^{4/5}\rho^{-1/10}
 +\gamma_{j+1}^{5/36}\rho^{-5/36}+r^{1/4}\rho^{-1/2}\right)F\nu^2\#\mathbb V.
 \end{aligned}
\]
On the other hand, Corollary~\ref{cor:coarse-frequency-final} gives
\begin{equation*}
    H(\rho)\lessapprox_{\varepsilon'}K^C\rho^{1/2}F\nu^2\#\mathbb V.
\end{equation*}
Using $\min\{a,b_1+b_2+b_3\}\leqslant\sum_{i=1}^3\min\{a,b_i\}$, we take geometric means of this bound with the three terms above in turn. The weights $1/6$, $5/23$, and $1/2$ on $\rho^{1/2}$, respectively, cancel the powers of $\rho$ and give
\begin{equation}\label{eq:layer-high-final}
 \begin{aligned}
 \frac{H(\rho)}{F\nu^2\#\mathbb V}
 &\lessapprox_{\varepsilon'} K^C
 \Biggl(
 \gamma_{j+1}^{2/3}\left(\frac{|Y|}{\gamma_j}\right)^{5/12}
 +\gamma_{j+1}^{5/46}\left(\frac{|Y|}{\gamma_j}\right)^{9/23}
 +\left(\frac{|Y|}{\gamma_j}\right)^{1/4}r^{1/8}
 \Biggr).
 \end{aligned}
\end{equation}

\textbf{Step 3.} Combining the bounds. Since $\gamma_j\lesssim1$, the first term in \eqref{eq:layer-high-final} is bounded by $(\gamma_{j+1}/\gamma_j)^{2/3}|Y|^{5/12}$. Combining \eqref{eq:layer-root-final},\eqref{4.20} and \eqref{eq:layer-high-final} and using $F\#\mathbb V\lesssim\#\mathbb W$ in the sum over $\rho,F,m$, we obtain
\[
 \begin{aligned}
 \frac{\|\Br_A Ef\|_{L^2(Y_j)}^2}{\nu^2\#\mathbb W}
 &\lessapprox_{\varepsilon'} K^C
 \Biggl(r^{13/36}
 +\left(\frac{\gamma_{j+1}}{\gamma_j}\right)^{4/5}|Y|^{5/12}\\
 &\hspace{5em}
 +\gamma_{j+1}^{5/46}\left(\frac{|Y|}{\gamma_j}\right)^{9/23}
 +\left(\frac{|Y|}{\gamma_j}\right)^{1/4}r^{1/8}\Biggr).
 \end{aligned}
\]
Combine the last two terms with \eqref{eq:layer-low-final} using weights $65/249$ and $5/21$ on $\gamma_{j+1}^{4/5}r^{2/5}$, respectively. This gives
\begin{align}
 \min\{\gamma_{j+1}^{4/5}r^{2/5},
        \gamma_{j+1}^{5/46}(|Y|/\gamma_j)^{9/23}\}
 &\leqslant
 \left(\frac{\gamma_{j+1}}{\gamma_j}\right)^{24/83}
 |Y|^{24/83}r^{26/249},\notag\\
 \min\{\gamma_{j+1}^{4/5}r^{2/5},
        (|Y|/\gamma_j)^{1/4}r^{1/8}\}
 &\leqslant
 \left(\frac{\gamma_{j+1}}{\gamma_j}\right)^{4/21}
 |Y|^{4/21}r^{4/21}.
 \label{eq:layer-remainders-final}
\end{align}
At $|Y|=r^{13/15}$, we have $|Y|^{5/12}=r^{13/36}$ and
\[
 |Y|^{24/83}r^{26/249}=r^{442/1245}\leqslant r^{13/36},
 \quad
 |Y|^{4/21}r^{4/21}=r^{16/45}\leqslant r^{13/36}.
\]
Both powers of $|Y|$ lie in $(0,5/12)$, so these expressions are bounded by $r^{13/36}$ below this threshold and by $|Y|^{5/12}$ above it. Using $\gamma_{j+1}/\gamma_j\leqslant r^\theta$, we have
\[
 \|\Br_A Ef\|_{L^2(Y_j)}^2
 \lessapprox_{\varepsilon'} K^C r^{O(\theta)}
 \left(r^{13/36}+|Y|^{5/12}\right)\nu^2\#\mathbb W .
\]
Choose $\theta$ and the power losses sufficiently small in terms of $\varepsilon$. Summing over the layers and the translated grids and applying \eqref{eq:packet-energy} gives \eqref{eq:floor-thirteen-mass}.
\end{proof}

\section{Proof of the main theorem}\label{sec:main-proof}

The broad operator uses the components indexed by the original $K^{-1}$-arcs. The following reduction is from \cite[Section~2.3, Method~2]{Wu25}.

\begin{lemma}[Two-ends reduction]\label{lem:wu-twoends}
Let $f=\sum_{T\in\mathbb W}f_T$ be a sum of $r$-scale wave packets. Let $X\subseteq B_r$ be a union of unit balls covered by a family $\mathcal Q$ of $r^{1/2}$-squares with $\|\Br_A Ef\|_{L^2(B)}\sim\lambda$ for each unit ball $B$ of $X$ and $|X\cap Q|\sim H$ for each $Q\in\mathcal Q$. Suppose $2\leqslant K\leqslant r^\kappa$, $A\geqslant r^{\varepsilon'}$, $0<\delta<1/10$, and $\varepsilon_0\ll\delta^2$. For sufficiently large $r$, one of the following holds with constants independent of $K$.
\begin{enumerate}
\item There are unions $X_k$ of unit balls from $X$ in balls of radius $r^{1-\delta}$, wave packet subsums $f_k$, and $A'\geqslant A/(\log r)^C$ such that
\begin{align}
 \|\Br_A Ef\|_{L^2(X)}^2
 &\lesssim(\log r)^C\sum_k
 \|\Br_{A'}Ef_k\|_{L^2(X_k)}^2, \label{eq:wu-concentrated}\\
 \sum_k\sum_{\sigma\in\Sigma}\|(f_k)_\sigma\|_2^2
 &\lesssim r^{C\delta^2}\sum_{T\in\mathbb W}\|f_T\|_2^2.
 \label{eq:wu-local-energy}
\end{align}
\item We have
\begin{equation}\label{eq:wu-twoends}
 \|\Br_A Ef\|_{L^2(X)}^2
 \lesssim r^{C\delta}|X|r^{-1/2}\sum_{T\in\mathbb W}\|f_T\|_2^2.
\end{equation}
\end{enumerate}
\end{lemma}

\begin{proof}[Proof of Proposition~\ref{main result}]
Fix $R$ and $K$. We prove by induction on $1\leqslant r\leqslant R$ that
\begin{equation}\label{eq:induction}
 \|\Br_A Ef\|_{L^2(X)}^2
 \leqslant C_{\eps,\eps'} R^\eps r^\eps
 |X|^{13/31}\sum_{\sigma\in\Sigma}\|f_\sigma\|_2^2,
 \quad A\geqslant r^{\eps'},
\end{equation}
for unions of unit balls $X\subseteq B_r$ with $r^{-1}X$ a Katz--Tao $(r^{-1},1)$-set. In the induction we allow $\supp f_\sigma$ to lie within distance $r^{-1/2}$ of $\sigma$.

By Lemma~\ref{pointwise estimate} and Cauchy--Schwarz,
\[
 \|\Br_A Ef\|_{L^2(X)}^2
 \lesssim |X|\sum_{\sigma\in\Sigma}\|f_\sigma\|_2^2.
\]
Since $|X|\lesssim r$, this proves \eqref{eq:induction} when $r\leqslant R^{\eps^2}$ and also covers bounded $r$ upon increasing $C_{\eps,\eps'}$. For $r>R^{\eps^2}$ we have $K\leqslant R^{\eps^4}\leqslant r^{\eps^2}$.

\textbf{Step 1.} Dyadic pigeonholing. Decompose each $f_\sigma$ into wave packets at scale $r$. Let $\mathbb T_r$ be the corresponding tube family. By dyadic pigeonholing and Lemma~\ref{Triangle inequality for broad operator}, choose dyadic numbers $\nu,\lambda,H>0$, a family $\mathbb W\subseteq\mathbb T_r$, a union $X_1\subseteq X$ of unit balls, a family $\mathcal Q$ of $r^{1/2}$-squares covering $X_1$ and $A_1\geqslant A/(\log r)^C$. Writing $f_1=\sum_{T\in\mathbb W}f_T$, we have
\begin{enumerate}
\item $\|f_T\|_2 \sim \nu$ for each $T\in\mathbb W$;

\item $\|\Br_{A_1}Ef_1\|_{L^2(B)} \sim \lambda$ for each unit ball $B$ of $X_1$;

\item $|X_1\cap Q| \sim H$ for each $Q\in\mathcal Q$.
\end{enumerate}
Moreover,
\begin{equation}\label{eq:initial-selection-induction}
    \|\Br_A Ef\|_{L^2(X)}^2
 \lessapprox \|\Br_{A_1}Ef_1\|_{L^2(X_1)}^2.
\end{equation}
By \eqref{eq:packet-energy} on each $\sigma$,
\begin{equation}\label{eq:selected-packet-energy}
 \sum_{\sigma\in\Sigma}\|(f_1)_\sigma\|_2^2
 \lesssim \sum_{T\in\mathbb W}\|f_T\|_2^2
 \lesssim \sum_{\sigma\in\Sigma}\|f_\sigma\|_2^2.
\end{equation}
Since $A_1\geqslant r^{\eps'/2}$ for large $r$, we apply Lemma~\ref{lem:wu-twoends} to $f_1$ and $X_1$ with $\delta=\eps^2$.

\textbf{Step 2.} One-end case. Let $A_2$ be the reduced parameter in \eqref{eq:wu-concentrated}. Then $A_2\geqslant r^{\eps'}/(\log r)^C \geqslant (r^{1-\eps^2})^{\eps'}$ for large $r$. Each $X_k$ inherits the Katz--Tao condition at radius $r^{1-\eps^2}$. The support condition also passes to $(f_k)_\sigma$ since $O(r^{-1/2})\leqslant (r^{1-\eps^2})^{-1/2}$ for large $r$. Applying \eqref{eq:induction} on the smaller balls and using \eqref{eq:wu-concentrated}--\eqref{eq:wu-local-energy} and \eqref{eq:initial-selection-induction}--\eqref{eq:selected-packet-energy} gives
\[
 \begin{aligned}
 \|\Br_A Ef\|_{L^2(X)}^2
 &\lesssim (\log r)^C C_{\eps,\eps'} R^\eps r^{(1-\eps^2)\eps}
 \sum_k |X_k|^{13/31}\sum_{\sigma\in\Sigma}\|(f_k)_\sigma\|_2^2\\
 &\lesssim r^{-\eps^3+O(\eps^4)} C_{\eps,\eps'} R^\eps r^\eps
 |X|^{13/31}\sum_{\sigma\in\Sigma}\|f_\sigma\|_2^2.
 \end{aligned}
\]
The factor $r^{-\eps^3+O(\eps^4)}$ closes the induction for small $\eps$ and large $r$.

\textbf{Step 3.} The two-ends case. Since $r^{-1/2}\ll K^{-1}$, each support $\supp(f_1)_\sigma$ meets $O(1)$ arcs of $\Sigma$. Split $\Sigma$ into $O(1)$ families $\Sigma_j$ so that each arc meets at most one such support from each family. Let $g_j=\sum_{\sigma\in\Sigma_j}(f_1)_\sigma$. By pigeonholing and disjointness within each family, there are $j_0$ and $A_1'\sim A_1$ such that
\[
 \|\Br_{A_1}Ef_1\|_{L^2(X_1)}^2
 \lesssim \|\Br_{A_1'}Eg_{j_0}\|_{L^2(X_1)}^2,
 \quad
 \sum_j\|g_j\|_2^2=\sum_{\sigma\in\Sigma}\|(f_1)_\sigma\|_2^2.
\]
Applying Proposition~\ref{prop:local-mass} to $g_{j_0}$ on $X_1$ with $A_1'\geqslant r^{\eps'/2}$ gives
\begin{equation}\label{2}
\begin{aligned}
 \|\Br_{A_1'}Eg_{j_0}\|_{L^2(X_1)}^2
 &\lesssim_{\eps,\eps'} K^C r^{\eps^2}
 \bigl(r^{13/36}+|X_1|^{5/12}\bigr)\|g_{j_0}\|_2^2\\
 &\lesssim_{\eps,\eps'} r^{O(\eps^2)}
 \bigl(r^{13/36}+|X_1|^{5/12}\bigr)
 \sum_{T\in\mathbb W}\|f_T\|_2^2.
\end{aligned}
\end{equation}
Combining \eqref{2} with \eqref{eq:wu-twoends}, \eqref{eq:initial-selection-induction}, and \eqref{eq:selected-packet-energy} gives
\[
 \begin{aligned}
 \|\Br_A Ef\|_{L^2(X)}^2
 &\lesssim_{\eps,\eps'} r^{O(\eps^2)}
 \min\left\{|X_1|r^{-1/2},\,r^{13/36}+|X_1|^{5/12}\right\}
 \sum_{\sigma\in\Sigma}\|f_\sigma\|_2^2.
 \end{aligned}
\]
The geometric mean with weights $13/31$ and $18/31$ gives
\[
 \begin{aligned}
 &\min\left\{|X_1|r^{-1/2},\,r^{13/36}+|X_1|^{5/12}\right\}\\
 &\quad\leqslant (|X_1|r^{-1/2})^{13/31}(r^{13/36})^{18/31}
 +|X_1|^{5/12}
 \lesssim|X_1|^{13/31}.
 \end{aligned}
\]
Since $|X_1|\leqslant |X|$ and $r^{O(\eps^2)}\leqslant r^\eps$ for small $\eps$, this completes the induction.
\end{proof}

\textbf{AI assistance.} GPT-6 Astra was used to assist with language editing, manuscript presentation, and the review of mathematical arguments.

\end{document}